\documentclass[11pt, oneside]{article}   	
\usepackage{geometry}                		
\usepackage{graphicx}				
\usepackage{amssymb}

\title{Some footnotes on Thurston's Notes \emph{The Geometry and Topology of 3-manifolds}}
\author{Athanase Papadopoulos\footnote{\emph{Institut de Recherche Mathématique Avancée} and
\emph{Centre de Recherche et d'Expérimentation sur l'Acte Artistique  (ITI CREAA)}. Address:
Université de Strasbourg and CNRS,
7 rue René Descartes,
67084 Strasbourg Cedex France;
email: athanase.papadopoulos@math.unistra.fr}}
 \date{\today}							

\begin{document}
\maketitle

\hfill \emph{J'aimais cette résistance coriace} 

\hfill \emph{dont je ne venais jamais à bout;}

\hfill \emph{mystifié, fourbu, je goûtais l'ambiguë}

\hfill \emph{volupté de comprendre sans comprendre :}

\hfill \emph{c'était l'épaisseur du monde}

\hfill  J.-P. Sartre, Les mots

\begin{abstract}

These are a few historical remarks, addenda and references with comments  on some topics discussed by Thurston in his notes \emph{The geometry and topology of three-manifolds}. The topics are mainly hyperbolic geometry, geometric structures, volumes of hyperbolic polyhedra  and the so-called Koebe--Andreev--Thurston theorem. I discuss in particular some works of Lobachevsky, Andreev and Milnor, with an excursus in Dante's cosmology, based on the insight of Pavel Florensky.

\bigskip

  The final version of this paper will appear as a chapter in the book \emph{In the tradition of Thurston, III} ed. K. Ohshika and A. Papadopoulos, Springer Nature Switzerland, 2024.

\bigskip

\noindent  {\bf Keywords:}
   geometric structure ; hyperbolic structure ;  $\mathcal{G}$-structure ; $(G,X)$-structure~;  hyperbolic volume ; volume of a hyperbolic tetrahedron ; volume of an ideal tetrahedron~; geometric convergence ; convex co-compact group ;  Koebe--Andreev--Thurston theorem ; Dante ; P. A. Florensky ; N. I. Lobachevsky.
 
 \bigskip
 
\noindent {\bf AMS codes:} 57-03 ;  57-06 ; 57K32 ;  57K35 ;  53C15  ; 01-160 ; 01A75 ; 01A20.

\end{abstract}

\medskip 

\vfill\eject

\tableofcontents

\section{Introduction}

 William Thurston's notes \emph{The geometry and topology of three-manifolds} recently appeared in print, as Volume IV of his Collected Works. 
The notes are associated with the courses he gave  at Princeton University during the academic year 1977-78 and the fall semester of 1978-79. Steve Kerckhoff writes in the
 Preface of this volume: ``Thurston started the first year by writing up notes himself, but after a while he enlisted Bill Floyd and me to help. [\ldots] What was written was a reasonable facsimile of the way Thurston was presenting his amazing mathematics to us, in real time."
The plan of the notes consists of 13 chapters, but Chapters 10 and 12 are missing, and Chapter 11 contains only a short section on extensions of vector fields.   The lectures corresponding to Chapter 7 were  given by John Milnor. The overall subject is geometric structures and the topology of three-manifolds. 

Thurston's notes require almost no prerequisites from the reader, they use only elementary mathematics, but they are difficult to read: one needs time, visualisation and imagination. On the other hand, reading every section and trying to understand the underlying ideas is always rewarding. It is not an exaggeration to say that hundreds of young mathematicians wrote their phD thesis by trying to understand some section of these notes.

   In this chapter, I shall make a few comments which I consider as footnotes to these notes, adding references here and there, sometimes mentioning some related work that preceded the notes,
  and,  more rarely, some developments and examples over the subsequent decades, needless to say, without proofs.
 
 My footnotes are divided into 6 sections, numbered 2 to 7.  
  Section \ref{s:NE}, which follows this introduction, consists of notes on the situation of non-Euclidean geometry at the time where Thurston's notes were written, and how these notes resurrected the subject. 
 Section \ref{s:Dante} starts with a sentence of Thurston: ```It would nonetheless be
distressing to live in elliptic space, since you would always be confronted with an image
of yourself, turned inside out, upside down and filling out the entire background
of your field of view."
This sentence reminds me of a passage in Dante's \emph{Divine Comedy}. Section \ref{s:Dante}  is then an excursion in the world of Dante, who imagined the universe as non-Euclidean, non-orientable and projective.
Section \ref{s:GS} is a note on the origin of Thurston's notion of geometric structure, highlighting the work of Ch. Ehresmann.
Section \ref{s:spheres} starts again with a sentence of Thurston on horospheres, and leads to an exposition of spheres and horospheres in hyperbolic space, in the work of N. I. Lobachevsky.
Section \ref{s:Andreev} consists of a few remarks on the work of E. M. Andreev to whom Thurston refers in his notes, regarding the classification of discrete isometry groups of hyperbolic space generated by reflections along hyperplanes with compact fundamental domain. This is the basis of  what became later the \emph{Koebe--Andreev--Thurston} theorem.
Section \ref{s:Volume} is concerned with the part of Thurston's notes that is due to  Milnor, on the computation of volume of hyperbolic polyhedra.

\section{Non-Euclidean geometry}\label{s:NE}

Chapter 2 of Thurston's Notes  is an introduction to the two non-Euclidean geometries of constant curvature, namely, spherical (or elliptic) and hyperbolic. 
Let me start by a few words on spherical geometry.

Spherical geometry is a classical subject, extensively studied in Greek Antiquity, in particular by Autolycus of Pitane (c. 360 -- c. 290 BC), Theodosius of Tripoli (c. 160 -- 90 BC) and Menelaus of Alexandria (c. 70 -- 140 AD).
 In these works, the sphere is considered as space having the same status as Euclidean 2-space, with lines, angles, etc. On the sphere, the lines are the great circles (that is, the intersection of the sphere with planes passing through the origin) and the angle between two lines is the dihedral angle between the planes that define them. The length of a segment joining two points is defined as the angle between the two Euclidean lines that join these points to the origin. With this, one defines triangles on the sphere and studies their properties, in particular trigonometry.  After  Greek Antiquity,  spherical geometry  was extensively studied by the medieval Arabs, who developed a complete trigonometric system and who discovered in particular duality and the notion of polar triangle. The modern period of spherical geometry starts with Euler and his collaborators and followers (Johann Lexell, Theodor von Schubert, Nikolaus  Fuss, etc.), who worked systematically on analogues in spherical geometry of theorems in Euclidean geometry due to Pappus and others. The references are \cite{Aujac-Autolycus, Menelaus, Pappus, Theodosius}.
 
Thurston writes in Chapter 2 of his notes: ``In the sphere, an object moving away from you appears smaller and smaller, until it reaches a distance of $\pi/2$. Then, it starts looking larger and larger and optically, it is in focus behind you. Finally, when it reaches a distance of $\pi$, it appears so large that it would seem to surround you entirely." 
 Rather than the sphere, Thurston considers the projective space, the quotient of the sphere by its canonical involution. Unlike the sphere, the projective space is a uniquely geodesic metric space: two points are joined by a single geodesic. Understanding the properties to which Thurston refers needs a mental effort.

After a glimpse into the mystery of elliptic geometry,\index{elliptic geometry} on which I will comment in the next section, 
  Thurston presents several models of hyperbolic geometry and he explains how to work with them, in particular  the sphere of imaginary radius, by which he means, the hyperboloid model.\footnote{The term ``sphere of imaginary radius" was known since the XIXth century, see later in this paragraph.), without reference to the hyperboloid model.} The importance of right angled hexagons is highlighted, and Thurston derives the basic trigonometric formulae for triangles, using this model.

One might remember that at the time Thurston wrote his notes, non-Euclidean geometry was a dormant subject; there were almost no textbooks on hyperbolic or spherical geometry, except for old books, mostly dating back to the nineteenth or the first years of the twentieth century. Let me mention here the interesting PhD thesis, defended in 1892, by L.  Gérard, a student of Paul Appell. The president of the jury of this thesis was Henri Poincaré. The dissertation is titled \emph{Sur la géométrie non euclidienne}. It was published in a book form \cite{Gerard}. In this book, the author derives the non-Euclidean trigonometric formulae by a model-free method, using essentially the fact that the angle sum of any triangle in the hyperbolic  plane is $<$ 2 right angles and the homogeneity of this plane. The same synthetic methods of deriving these formulae work for spherical trigonometry. I also mention  Barbarin's \emph{La géométrie non euclidienne} \cite{Barbarin}, a booklet in which the author also establishes by a model-free method the trigonometric formulae. This booklet  contains a section on neutral geometry, also called absolute geometry, that is, a geometry whose axioms consist of those of Euclidean geometry with the parallel axiom deleted. This is the geometry whose propositions are valid both in Euclidean and hyperbolic geometry. In fact, it was already noted since Greek Antiquity that a large number of statements of Euclidean geometry hold without the parallel axiom, that is, in absolute geometry, and it was even thought that the parallel axiom was useless, namely, that it follows from the other axioms.  This is the main reason for the large effort that was spent since those days to prove that the parallel axiom is a consequence of the others. Some of the greatest mathematicians of the modern period, including Euler, Lagrange, Legendre, Fourier and others, were also involved in this effort; they all thought  that the parallel axiom was a consequence of the others, and therefore, that there is no hyperbolic geometry. On this modern period, see e.g. the report in \cite[p. 24-28]{Papa-Lambert}.
Later books include Bonola’s \emph{Non-Euclidean geometry, a critical and historical study of its developments} \cite{Bonola}. One might add the book \cite{FN}  by Werner Fenchel and Jakob Nielsen, which appeared in print in 2003 but which was circulated before, in the form of mimeographed notes, under the name “the Fenchel--Nielsen manuscript”. These notes developed through lectures given by the authors at the University of Copenhagen for more than 40 years. A first draft was already planned to appear in the Princeton Mathematical Series in 1948, but then, as the subject grew up, the authors made several revisions, and the book never appeared before Nieslen’s death (1959). Fenchel continued this writing project alone until his death (1988). The book was published in 2003, edited by Asmus Schmidt.

Regarding the classical textbooks on hyperbolic geometry, I also mention
Hadamard's \emph{Leçons de géométrie élémentaire} \cite{Hadamard}, (in principle a textbook for last-year high-school students), which contains a chapter on hyperbolic geometry. Finally, in some books on complex analysis,  e.g. in Carathéodory's \emph{Theory of functions of a complex variable} \cite{Cara}, there are also chapters on spherical and hyperbolic geometry. 
 
When Thurston wrote his notes, hyperbolic geometry was very rarely (if ever at all) taught in universities except in the Soviet Union, Romania and some other Eastern European countries.  I personally heard of it at a course on complex analysis by V. Poénaru  (who was trained in Romania), when I was an undergraduate at Orsay. I found at that time, in a Russian bookstore in Paris, an excellent textbook by N. Efimov, titled  \emph{Géométrie supérieure} \cite{Efimov},  which contains a large section on hyperbolic geometry.  
 In any case, I think that for most of the first readers of  Thurston's notes, hyperbolic geometry was a new subject.   Spherical geometry  was also poorly studied in the twentieth century.
 
 Much later, I  got acquainted with Lobachevsky's\index{Lobachevsky, Nikolai Ivanovich} works. The latter used two different words for  hyperbolic geometry:  \emph{Imaginary geometry}  \cite{Lobachevsky-Imaginary, Lobachevsky-Imaginaire} and \emph{Pangeometry} \cite{Pangeometry}. The word \emph{imaginary} refers to the sphere of imaginary radius. The reason Lobachevsky\index{Lobachevsky, Nikolai Ivanovich}  talked about a sphere of imaginary radius is that after he derived the trigonometric formulae, he realized that these formulae are the same as those of spherical geometry provided one replaces the terms involving the cosine of length of a segment by its hyperbolic cosine. Since at the formal level this consists in replacing terms like $\cos x= \frac{e^{ix}+e^{-{ix}}}{2}$ by  $\cos x= \frac{e^x+e^{-x}}{2}$, Lobachevsky concluded that in hyperbolic geometry one works on a sphere of imaginary radius. He also used the name \emph{Pangeometry}\index{Lobachevsky, Nikolai Ivanovich!\emph{Pangeometry}}  because he saw that hyperbolic geometry contains the two other geometries: the spherical, and the Euclidean. Indeed, he showed that  a geometric sphere in a three-dimensional hyperbolic space carries the same geometry of the usual sphere---the model of spherical geometry---and that the Euclidean plane sits in the three-dimensional hyperbolic space as a horosphere. I will discuss this in \S \ref{s:spheres}.
 
   Lobachevsky's\index{Lobachevsky, Nikolai Ivanovich}  writings contain most of what we need to know in hyperbolic geometry, even though there were no Euclidean models available to him. He derived a complete system of hyperbolic trigonometric formulae, and he asserted at several occasions that once we know the trigonometric formulae of some geometry, then we have all the tools for understanding that geometry. Indeed,  from the trigonometric formulae, we have a complete knowledge of the geometry of triangles, and the geometry of the triangles gives complete information on the geometry of the space, both infinitesimally and in the large. 
   
   Thurston was knowledgeable in hyperbolic geometry since early in his career. In his 1972 paper on foliations \cite{Thurston1972}, the hyperbolic plane already appears as an essential tool. Let me also quote a passage from Sullivan's article in \cite{LP}, which shows that a little bit before 1970, Thurston was already thinking about hyperbolic surfaces obtained by gluing geodesic hyperbolic octogons along their boundaries. The constructions of the universal cover of such a surface seemed completely new for others in Berkeley. Sullivan writes:
   \begin{quote}\small 
   This happens a little bit before 1970:   [\ldots]  I  had spoken to my old
friend Mo Hirsch about  Bill Thurston who was working with Mo and was
finishing in his fifth year after an apparently slow start.  Mo or someone
else  told how  Bill's oral exam was a slight problem because when asked
for an example of the universal cover of a space Bill chose the surface of
genus two and started drawing awkward octogons with many edges coming
together at each vertex. This exposition quickly became an unconvincing
mess on the blackboard. I think Bill was the only one in the room who had
thought about this  nontrivial kind of universal cover.
\end{quote}

At Orsay, Thurston's notes used to arrive in installments, soon after they were written. Norbert A'Campo, in \cite{LP}, mentions that, besides Poénaru, only Douady and Deligne knew about this subject. He recalls that Siebenmann asked him to give a course on hyperbolic geometry, and that his response was that he knew very little about it. It was a happy coincidence that A'Campo visited the Mittag-Leffler Institute, just after Siebenmann's request, and there, in one of the attics  he 
found a set of old papers, notes and documents on hyperbolic geometry.  He went through them and took some notes, and he came back to France with enough material to build a course. This is how, after Thurston's notes arrived, a course on hyperbolic geometry was taught at Orsay. The full story is recounted in \cite{Papa-Cafe}. Several years later, I went with A'Campo through his notes, we expanded them together, and we published the result, see \cite{ACP}.
A'Campo writes in \cite{LP} that he asked Thurston how he came to know about hyperbolic geometry,
and Thurston's answer was that he first learned it from his father.

Finally, let me note that Lobachevsky,\index{Lobachevsky, Nikolai Ivanovich}  in his courses at Kazan University and in his first memoirs, considered at the same time spherical and hyperbolic geometry; see \cite{Lobachevsky-New}.

Needless to say, several textbooks on non-Euclidean geometry were written by various authors after Thurston's notes appeared.

 \section{Excursus: Dante}\label{s:Dante}

This section is motivated by a sentence in Chapter II of Thurston's Notes, in which the latter introduces the two non-Euclidean geometries: the hyperbolic and the spherical. Of the latter, he considers the non-orientable version, that is, elliptic geometry. 
He writes:  ``It would nonetheless be
distressing to live in elliptic space, since you would always be confronted with an image
of yourself, turned inside out, upside down and filling out the entire background
of your field of view." My excursion concerns the history of this sentence. Most of what I write here about this question is extracted from my article \cite{Papa-Florensky}.

In 1922, the Russian mathematician, philosopher and theologian Pavel Florensky\footnote{Pavel Aleksandrovich Florensky (1882-1937)\index{Florensky, Pavel Aleksandrovich} was a mathematician, physicist, chemist, botanist, theologian, engineer, philosopher, art historian, specialist in electrodynamics, marine botany, linguistics,  languages, symbol, esoterism, gnoseology, logic, aesthetics, semiotics and  other fields of knowledge. He was among the Russian intellectuals of the Stalin period who spent the last years of their life in the Gulag. Those who knew his work called him the Russian Leonardo, or the Russian Pascal. Solzhenitsyn wrote that Florensky was ``perhaps the most remarkable person devoured by the Gulag". He is the author of an immense body of works, which is now starting to be published in various Western languages. Cédric Villani, in the preface to the French edition of his \emph{Imaginaries in geometry}\index{Florensky, Pavel Aleksandrovich!\emph{Imaginaries in geometry}} \cite{Imaginaires}, described him there as ``one of the most original thinkers of his times".}
published a booklet
titled \emph{Imaginaries in geometry}. The booklet appeared  recently in French translation, with a Preface by Cédric Villani. 
 The work is partly on complex numbers (as the title suggests), partly  on imagination, in the sense of Hilbert--Cohn-Vossen and of Thurston's \emph{Geometry and imagination} \cite{Hilbert-Imagination, Thurston-Imagination}, and partly on imagination as common people intend it. I have commented on this booklet and on other works of Florensky in the article \cite{Papa-Florensky}.
 In \S 9 of this booklet, Florensky mentions a passage from Dante's \emph{Divine Comedy}\index{Dante!\emph{Divine Comedy}} in which the Italian poet describes 
   an imaginary journey he undertakes, with the objective of meeting the Creator. The journey takes him successively to Hell, Purgatory and Paradise. This is not the place to comment on the poem, but let me mention only that this journey, which has to be interpreted as a fable, epitomizes the Poet's exile (he had been expelled from Florence by the city's authorities), and it also recounts the journey of each and every one of us. Readers interested in the history of ancient science can follow the Poet's description of his journey, during which he presents the Ptolemaic vision of the universe. 
   
In the last sixty years, several mathematicians and physicists commented on a passage of the  \emph{Divine Comedy}, sometimes at length,
stating that Dante describes there the world as a three-dimensional sphere. They did so forgetting an  important observation made by Florensky, the one that made me think of Thurston's sentence I mentioned above. I shall comment on this observation by Florensky. Let me first mention some of the  mathematicians who have talked about Dante's description of the topology of the universe. The  complete list is long, and I shall give just part of it.
I will mention these authors approximately in chronological order: 
 The physicists I. Ozsv\'ath, and E. Schücking,  in the article  \emph{An anti-Mach metric}, published in the book \emph{Recent developments in general relativity}  \cite{O-Sch} (1962), and the same authors 34 years later, in 
the article \emph{The world viewed from outside}, published in the Journal of Geometry and Physics \cite{O-Sch-Dante}  (1998); 
  M. Peterson, in his paper \emph{Dante and the 3-sphere}, published in the American Journal of Physics \cite{Peterson-Physics} (1979),  
  and the same author some  30 years later, in his paper  \emph{The geometry of paradise} published in the Mathematical Intelligencer \cite{Paterson-Intelligencer}  (2008); 
   R. Osserman, in his book  \emph{Poetry of the Universe: A Mathematical Exploration of the Cosmos} \cite{Osserman} 
  (1995) and in his article
\emph{Curved Space and Poetry of the Universe}  published in the book \emph{The book
of the cosmos : imagining the universe from Heraclitus to Hawking} \cite  {Osserman2000}  (2000); 
   W.  Egginton, in his paper  \emph{On Dante, hyperspheres, and the curvature of the medieval cosmos}  published in the Journal of the History of Ideas \cite{Egginton}  (1999);
M. Wertheim, in her book \emph{The pearly gates of cyberspace: A history of space from Dante to the internet} \cite{Wertheim} (1999);  G. Mazzotta,  in his article  \emph{Cosmology and the kiss of creation (Paradiso 27-29)} published in the  Dante Studies \cite{Mazzotta}  (2005); 
  S. L. Lipscomb in his book  \emph{Art meets mathematics in the fourth Dimension} \cite{Lipscomb}  (2011); 
Curt McMullen in his survey, \emph{The evolution of geometric structures on 3-manifolds} \cite{McMullen-Dante}, published in the Bulletin of the American Mathematical Society (2011);\footnote{McMullen's survey starts  with the sentence:
 ``In 1300, Dante\index{Dante!cosmology}  described a universe in which the concentric terraces of hell-nesting down to the center of the earth-are mirrored by concentric celestial spheres,
rising and converging to a single luminous point. Topologically, this
finite yet unbounded space would today be described as a three-dimensional sphere."}
 M. Bersanelli, in his article \emph{From Dante's\index{Dante!\emph{Divine Comedy}}  universe to contemporary cosmology}, published in Rend. Scienze Ins. Lombardo \cite{Bersanelli}  (2016);
C. Rovelli, in his article \emph{Michelangelo's stone: an argument against Platonism in mathematics}, published in the European Journal for Philosophy of Science
\cite{Rovelli} (2017);
E. Ghys, in  his book \emph{A singular mathematical promenade} \cite{Ghys} (2017);
K. Sunada, in the article \emph{From Euclid to Riemann and beyond: how to describe the shape of the universe}, published in a book I edited  \cite{Sunada}  (2019). 
Let me finally mention  the  interesting and beautifully illustrated article \cite{Dreyer} by J. L. E. Dreyer, \emph{The cosmology of Dante}, published\index{Dante!cosmology}  in 1921 in \emph{Nature},  in which the author surveys Dante's cosmological ideas expressed  in his poem (with no mention of the topology of the space described there).
The article ends with the words: 
``To the student of the history of science it is a never-failing source of pleasure to find medieval cosmology so beautifully illuminated in the writings of the great Florentine poet" \cite{Dreyer}.

I consider that most of the authors I just quoted failed in two ways:
First of all, none of them refers to the work of Florensky\index{Florensky, Pavel Aleksandrovich!\emph{Imaginaries in geometry}} on this subject, which was published between 40 years and almost a century before.\footnote{First edition, Moscow, Pomor'ye 1922 (1000  printed copies), new Russian edition, Moscow, Lazur, 1991.} For several among them, this is understandable because Florensky's booklet, until its first translation in French in 2016, existed only in Russian. Let me mention incidentally that one rarely finds authors on this list who quote others on the same list who made equivalent statements. Among
the few authors who cite others are Rovelli, Sunada and Ghys, who quote Peterson.

The second way I think these authors failed is that all of them, with the goal of highlighting the universe's topology described by Dante in his poem, declared that the space which is described in the \emph{Divine Comedy}  is the 3-dimensional sphere, but they missed the part on non-orientability. This is the property highlighted by Florensky.
The only paper I am aware of where the authors mention non-orientability is a report on Florensky's\index{Florensky, Pavel Aleksandrovich!\emph{Imaginaries in geometry}} commentary on Dante's journey, by Bayuk and Ford, titled \emph{Dante's cosmology revisited} and published in the  Archives Internationales d'Histoire des Sciences \cite{BF}.  

  Section 9 of the \emph{Imaginaries},\index{Florensky, Pavel Aleksandrovich!\emph{Imaginaries in geometry}} which is the last section of the booklet, is dedicated to the 600th anniversary of the death of Dante. This event was celebrated in Russia on September 14, 1921.  I will comment now on this section of Florensky's work.  Let me first recall Dante's setting.

  \begin{figure}[htbp]
\centering
\includegraphics[width=14cm]{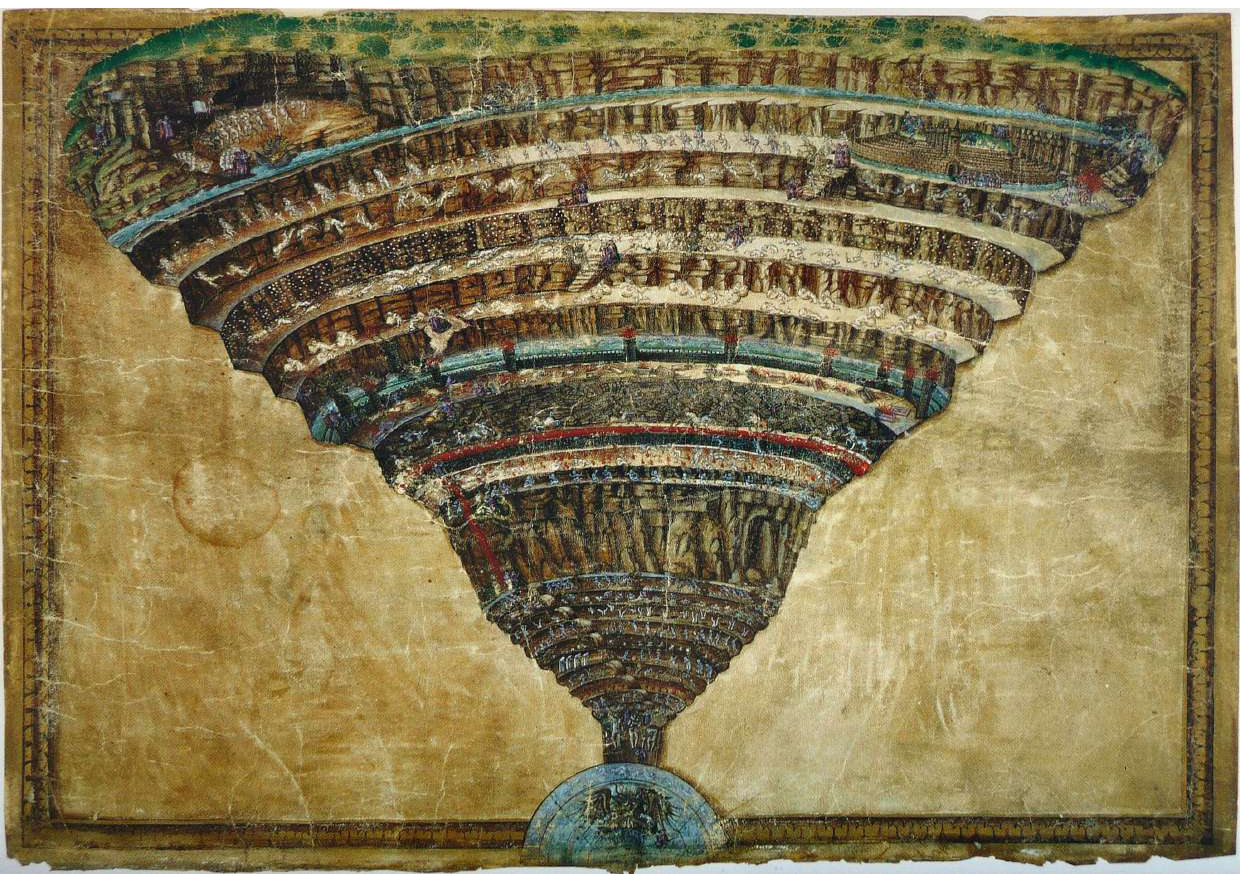}
\caption{A map of the inferno, funnel-shaped and made out of nine circles of torment; from the circle of largest to the smaller  diameter, they are the realm of  limbo, lust, gluttony, greed, anger, heresy, violence, and the last two circles are for fraud: flatterers and treaters. Lucifer dwells in the smallest (and deepest) circle. From the Divine Comedy Illustrated by Sandro Botticelli  (1445-1510),  between 1480 and 1490. Vatican Library.}
\label{fig:Botticelli}
\end{figure}

\begin{figure}[htbp]
\centering
\includegraphics[width=8cm]{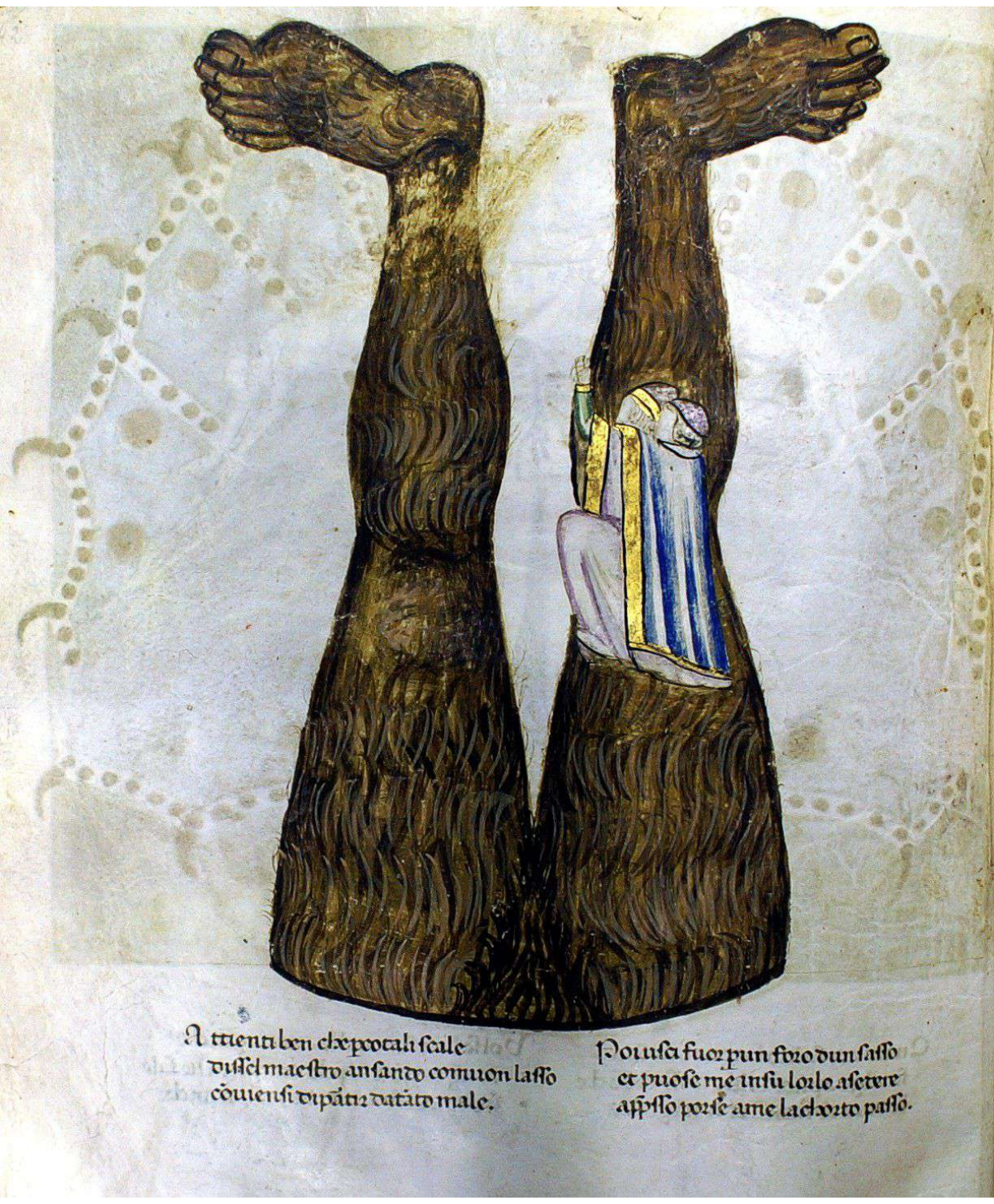}
\caption{Lucifer upside down. Codex Altonensis Dante, Commedia, 14th c., Bibliotheca Gimnasii Altonani, Hamburg.}
\label{fig:Lucifer}
\end{figure}

 Dante is accompanied in his journey by the poet Virgil, the Roman poet from the 1${}^{\mathrm{st}}$ century {\sc bc}, author of the \emph{Aeneid}, who guides him through Hell and most of the Purgatory.   
 The two poets start at Florence, for a long descent to the bottom of the Earth, taking the funnel-shaped slopes (see Figure \ref{fig:Botticelli})  leading to the narrowest circle of Hell inhabited by Lucifer, the angel who was cast down from heaven and who became the Master of the kingdom of darkness. During all the time of their descent, the heads of the two poets are directed towards their place of departure, that is, upwards, and their feet are directed towards the center of the Earth. But when they find themselves near Lucifer, they see him with his head down. It is also conceivable that Lucifer is standing up on his feet, and that they are themselves turned upside down, with their heads towards the bottom of the Earth and their feet towards the opposite side. We have reproduced in Figure \ref{fig:Lucifer} a   14${}^{\mathrm{th}}$ century representation of Lucifer and one of the 
the two poets, upside down with respect to each other.

Now comes geometry. 

Florensky,\index{Florensky, Pavel Aleksandrovich} refers 
to a passage of the \emph{Divine Comedy}\index{Dante!\emph{Divine Comedy}}  where the Florentine poet speaks of triangles whose  angle sum is not equal to two right angles.  In this passage, Dante\index{Dante!\emph{Divine Comedy}}  quotes King Solomon who asks the Lord whether there exists a triangle inscribed in a semi-circle, with one side being the diameter, and no right angle  (\emph{Paradiso},  Canto XIII, v. 101-102):\footnote{transl. H. W. Longfellow \cite{Dante}.}
\begin{verse}
non \emph{si est dare primum motum esse}, \hfill
nor \emph{si est dare primum motum esse},\footnote{``nor if we have to admit a prime motion" (Dante wrote this line in Latin).}\\

o se del mezzo cerchio far si puote \hfill  	nor if in a semicircle a triangle can be formed \\

triangol s\`\i \  ch'un retto non avesse.\hfill   without its having one right angle. \\
 
\end{verse}

Florensky makes the comment that ``the image of this universe cannot be represented by Euclid's geometry, just as Dante's metaphysics is incommensurable with Kant's philosophy "\cite[p. 71]{Imaginaires}).\footnote{My translation from the French.} He\index{Florensky, Pavel Aleksandrovich} notes that three authors before him made the observation that Dante's world is non-Euclidean, namely,  George Bruce
Halstead (1853-1922),   Heinrich Martin Weber (1842-1913), and   Maximilian Simon
(1844-1918).

 Dante's\index{Dante!cosmology}  discussion on non-orientability needs a longer quote and is also based on Florensky's analysis, which I will review now.

Florensky\index{Florensky, Pavel Aleksandrovich} starts by recalling Dante's description of the Ptolemaic model of the universe, composed of interlocking spheres: the Earth, the celestial bodies and the spheres on which they lie, together with the Empyrean stratum, a sphere surrounding the  whole universe and occupied by the element of fire. 
Here is the passage of the\index{Dante!\emph{Divine Comedy}}  \emph{Divine Comedy} referred to  by Florensky\index{Florensky, Pavel Aleksandrovich}  (\emph{Inferno}, Canto XXXIV, v. 76-96):

\begin{verse} \small\small
Quando noi fummo là dove la coscia  \hfill When we were come to where the thigh  revolves \\

si volge, a punto in sul grosso de l'anche,  \hfill  Exactly on the thickness of the haunch,\\

lo duca, con fatica e con angoscia,    \hfill  The Guide, with labour and with   hard-drawn  breath, \\   

\medskip                        

volse la testa ov'elli avea le zanche,  \hfill Turned round his head where he had had his legs, \\

e aggrappossi al pel com'om che sale, \hfill    And grappled to the hair, as one who mounts, \\

s\`\i \  che 'n inferno i' credea tornar anche.   \hfill    So that to Hell I thought we were  returning. \\    

\medskip                  

``Attienti ben, ché per cotali scale",\hfill    ``Keep fast thy hold, for by such stairs as   these," \\

disse 'l maestro, ansando com'uom lasso, \hfill   The Master said, panting as one  fatigued, \\

``conviensi dipartir da tanto male".  \hfill  ``Must we perforce depart from so much  evil." \\   

\medskip                              

Poi usc\`\i \  fuor per lo f\'oro d'un sasso, \hfill  Then through the opening of a rock he  issued,  \\     

e puose me in su l'orlo a sedere;  \hfill  And down upon the margin seated me; \\     

appresso porse a me l'accorto passo.  \hfill     Then tow'rds me he outstretched his  wary  step.  \\      

\medskip

Io levai li occhi e credetti vedere \hfill   I lifted up mine eyes and thought to see  \\     

Lucifero com'io l'avea lasciato,\hfill     Lucifer in the same way I had left him;  \\     

e vidili le gambe in sù tenere;       \hfill       And I beheld him upward hold his  legs.     \\

\medskip

e s'io divenni allora travagliato,\hfill   And if I then became disquieted,  \\     

la gente grossa il pensi, che non vede \hfill   Let stolid people think who do not see  \\     

qual è quel punto ch'io avea passato.   \hfill      What the point is beyond which I had passed.      \\       

\medskip

``Lèvati sù", disse 'l maestro, ``in piede: \hfill  ``Rise up," the Master said, ``upon thy  feet;  \\     

la via è lunga e 'l cammino è malvagio,  \hfill The way is long, and difficult the road,  \\     

e già il sole a mezza terza riede".  \hfill   And now the sun to middle-tierce returns."

\end{verse}

    \begin{figure}[htbp]
\centering
\includegraphics[width=8cm]{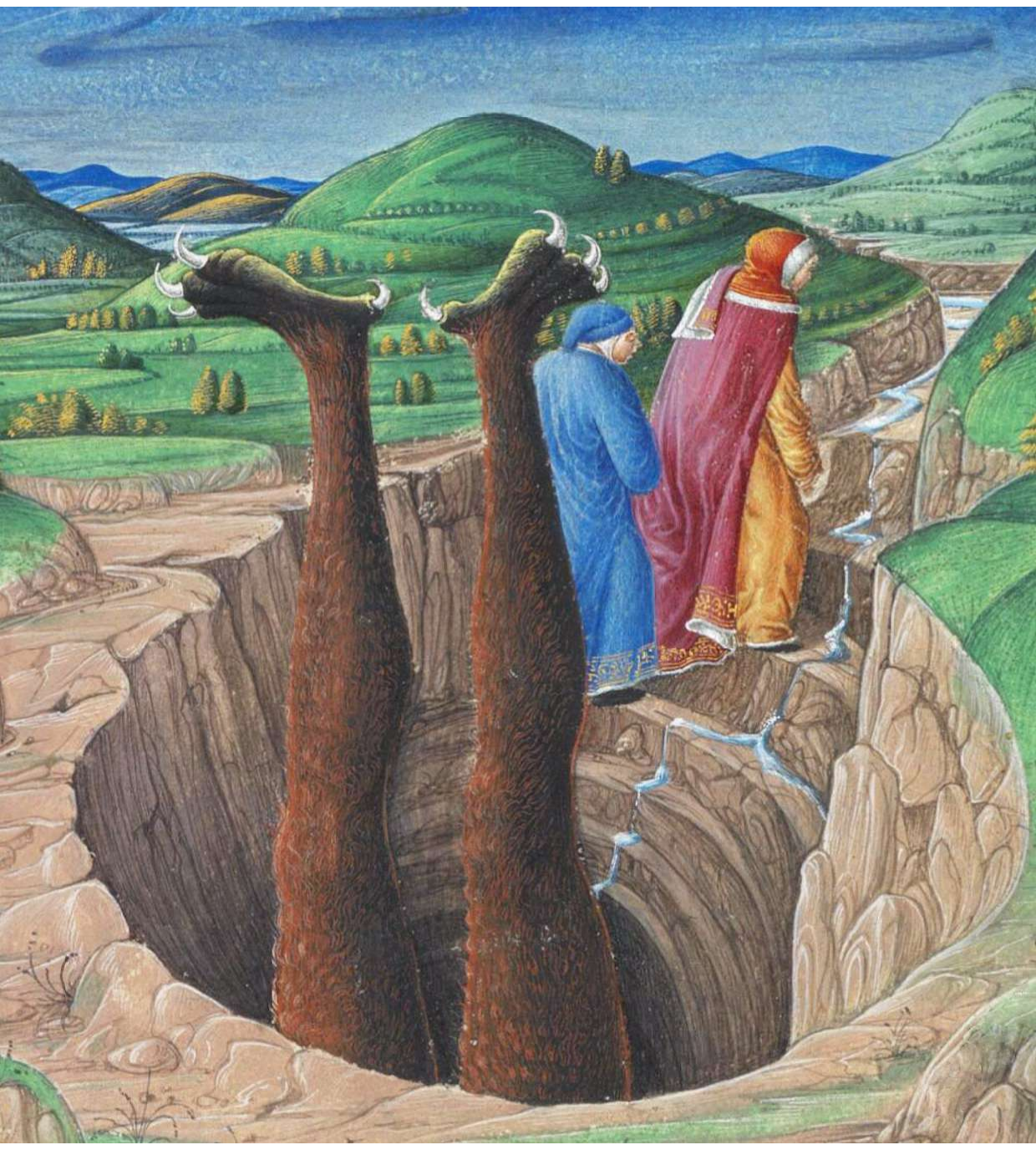}
\caption{Illustration of canto XXXIV of Dante's \emph{Inferno}, with Lucifer upside down, and the two poets as they return to the clear world. From a manuscript in Rome, Vatican Library, Ms Urb. Lat. 365, fol. 95v. Ferrara, 1474-1482.}
\label{fig:Ferrara}
\end{figure}

%
%

The two poets, after having avoided the center of the universe, continue, following the same path, but now they find themselves moving upwards, until they traverse a crater track, to end up in the clear world. This is illustrated in Figure \ref{fig:Ferrara}, reproduced from a manuscript in the Vatican Library. In Dante's description (Inferno, Canto XXXIV, v. 130-139): 

\begin{verse}\small
 
Lo duca e io per quel cammino ascoso \hfill The Guide and I into that hidden road \\

intrammo a ritornar nel chiaro mondo;    \hfill   Now entered, to return to the bright world; \\

e sanza cura aver d'alcun riposo,     \hfill     And without care of having any rest   \\    

\medskip                        

salimmo sù, el primo e io secondo,    \hfill  We mounted up, he first and I the second,  \\

tanto ch'i' vidi de le cose belle     \hfill  Till I beheld through a round aperture \\

che porta 'l ciel, per un pertugio tondo.     \hfill  Some of the beauteous things that Heaven 
 doth bear; \\

\medskip

E quindi uscimmo a riveder le stelle.   \hfill Thence we came forth to rebehold the  stars.  \hfill

\end{verse}

Dante\index{Dante!\emph{Divine Comedy}} begins his ascent towards the celestial spheres, always following the same line.
After a passage through the Empyrean where he can contemplate the Glory of God, who dwells in the last sphere after the Empyrean, he  finds himself  again in Florence, where he began his journey.

\begin{figure}
\begin{center} 
\includegraphics[width=10 cm]{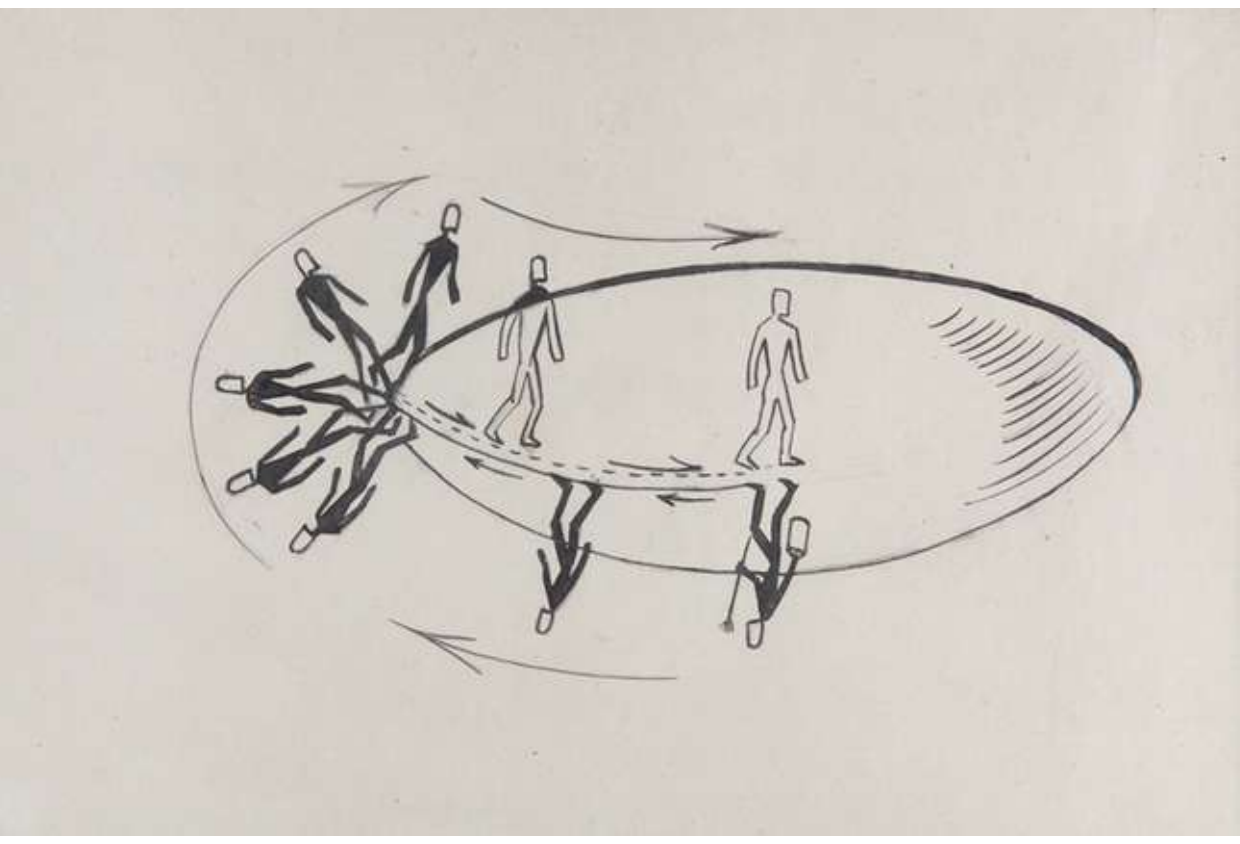} 
\caption{Illustration from Florensky's \emph{Imaginaries in geometry}.} 1922. Paper, ink
\label{fig:Imaginary}
\end{center}
\end{figure}

Thus, Dante, having traveled along a (non-Euclidean) straight line in the 3-dimensional world, finds himself at the same place and in the same position, but after having undergone a complete reversal. He has followed a geodesic, in this non-Euclidean space, and he has experienced an orientation-reversal. Obviously, Dante did not have the technical means to express this in the way a mathematician would today.
 Figure \ref{fig:Imaginary} is reproduced from Florensky's \emph{Imaginaries in geometry}; it represents a man walking in a non-orientable world. 
Florensky continues \cite[\S 9, p.73]{Imaginaires}, concerning Dante's journey:

\begin{quote}\small

Thus, after traveling all the way along a straight line and turning himself around once during his journey, the poet returns to the place where he was before, in the same position in which he left. 
Consequently, had he not turned around during his journey, he would have returned in a straight line to his starting point, head down and feet up. This means that the surface on which Dante moves is such that a straight-line movement with a single change of direction returns the body to the starting point in a straight position, and a straight-line movement with no change of direction returns the body to the starting point in an inverted position. 
Evidently, this surface is 
a {\bf Riemann surface}\footnote{The bold face letters are Florensky's.}   because the lines it contains are closed,\footnote{In the nineteenth century and at the beginning of the twentieth, mathematicians used to call the sphere, with its round metric a \emph{Riemann surface} or a \emph{Riemannian space}. Today, the terminology \emph{Riemann surface}  is used for a sphere (or another surface) equipped with a conformal structure.} and (2) one-sided because it turns the perpendicular around as it moves along it.\footnote{Presumably, a one-sided Riemann surface is a projective plane.} These two circumstances suffice for the geometric characterization of  {\bf Dante's space
as being built according to  the type of elliptic geometry}. One must recall that {\bf Riemann},
using methods of differential geometry, was unable to examine the {\bf global} form of surfaces.
This is why the object of his geometrical discussions was indifferently two geometries that are far from identical, one based on the elliptic plane, and the other on the spherical plane. In 1871, Felix Klein noted that the spherical plane has the character of a bilateral surface, whereas the elliptical plane has that of a unilateral surface. Dante's space closely resembles elliptical space. In this way, the medieval concept of the finiteness of the world is revealed in an unexpected way. 
\end{quote}

Let me conclude this section by noting that Florensk's \emph{Imaginaries} \cite{Imaginaires}, after the passage on Dante, contains a discussion of  Einstein's\index{relativity theory}  theory of relativity, and that this is an important component of the booklet.   There is a lot to say about this, and I refer the interested reader to the edition \cite{Imaginaires} and to the commentaries contained in that booklet. Let me mention at least that Dante's cosmology, in relation with Einstein's theory of relativity, is discussed in a paper by  J. Callahan,  titled \emph{The curvature of space in a finite universe} \cite{Callahan}, that appeared in the Scientific American referring to the book \emph{Klassische Stücke der Mathematik} \cite{Speiser}. Callahan writes: ``The spiritual
world completes the material world exactly as one viewing screen of Einstein's
model of the galactic system completes the other. The overlap between the two
is revealed by the correspondence between the celestial spheres and their angelic counterparts, and again as in Einstein's model the farther a sphere is from
the center of one chart, the nearer its
counterpart is to the center of the other.
And the speeds of the material spheres
and of the spiritual spheres are in harmony."

\section{Geometric structures: From Ehresmann and Haefliger  to Thurston} \label{s:GS}

 Chapter 3 of Thurston's notes  is concerned with geometric structures on manifolds.  In this chapter, Thurston introduces two slightly different notions of geometric structures:  $\mathcal{G}$-manifolds and $(G, X)$-manifolds, the second one being a special case of the first.

The definition uses the notion of pseudo-group. One starts with a collection $\mathcal{G}$ of local homeomorphisms of some manifold satisfying four natural properties making $\mathcal{G}$ a pseudo-group: 

\begin{enumerate}

\item Restriction: If $g\in \mathcal{G}$, then $g$ restricted to any open subset of its domain is also in $\mathcal{G}$.

\item  Composition: If $g_1, g_2\in \mathcal{G}$, then $g_1\circ g_2$, on its domain of definition, is also in $\mathcal{G}$.

\item Inverse: If $g\in \mathcal{G}$,  then $g^{-1}$ is also in $\mathcal{G}$.

\item Being in $\mathcal{G}$ is a local property: If  $U=\cup_{\alpha} U_{\alpha}$ and if $g:U\to V$ is a local homeomorphism whose restriction to each $\alpha$ is in $\mathcal{G}$, then $g\in \mathcal{G}$.

\end{enumerate}

A  $\mathcal{G}$-manifold\index{Gmanifold@$\mathcal{G}$-manifold} (or a $\mathcal{G}$-structure\index{Gstructure@$\mathcal{G}$-structure} on a manifold) is then a manifold obtained by piecing together in a coherent way elements of $\mathcal{G}$. A particular case of a $\mathcal{G}$-structure\index{Gstructure@$\mathcal{G}$-structure} is the case where one starts with a manifold $X$ equipped with the action of a group $G$, and takes $\mathcal{G}$ to be the pseudo-group of local homeomorphisms between open subsets of $X$ determined by this action. The resulting $\mathcal{G}$-manifold is what Thurston calls a $(G, X)$-manifold.

For example if $X=\mathbb{R}^n$ and $\mathcal{G}$ the pseudo-group of restrictions of local $C^r$ for some $r\geq 1$ (or piecewise linear, etc.) homeomorphisms  of $\mathbb{R}^n$, then the resulting $\mathcal{G}$-manifold is (by definition) a  $C^r$- (or \emph{piecewise linear}, etc.)-\emph{manifold}. 
 Another typical example of a geometric structure is a locally homogeneous  complete  Riemannian metric on a manifold, or, more generally, a locally homogeneous complete metric space. Here, a metric space is called locally homogeneous if any two points have neighborhoods that are isometric by an isometry sending one of these points to the other one.

Interestingly, Thurston refers to a paper by\index{Haefliger, André} Haefliger \cite{Haefliger1971} for the fact that in the pseudogroup definition of a $\mathcal{G}$-manifold that we recalled above, one may restrict the definition to a pseudogroup acting on pieces of $\mathbb{R}^n$, instead of acting on a general manifold, provided $\mathcal{G}$ is transitive. Haefliger's paper concerns the case of codimension $q$ foliations on an $n$-dimensional manifold. In particular, Haefliger addresses the question of finding 
  conditions under which a field of codimension $q$ tangent planes in a manifold is homotopic to a field tangent to a foliation.\footnote{The question was already asked by G. Reeb in his paper \cite{Reeb1952} (1952).} It may be useful here, for the sake of understanding the history, to remember that Thurston, before starting his work on the geometry and topology of three-manifolds, used related notions of geometric structures in his work on foliations. In fact, Haefliger's papers are already quoted in Thurston's PhD dissertation \cite{Thurston-dissertation} (1967). Haefliger, in his  1957 paper \cite{Haefliger1970} titled \emph{Feuilletages sur les variétés ouvertes } and which is already quoted by Thurston in his PhD dissertation, introduces the notion of \emph{$\Gamma_q$-structure}, where $\Gamma_q$ is the pseudo-group of local $C^{\infty}$-diffeomorphisms of $\mathbb{R}^q$, which  bears a certain analogy to that of geometric structure. Haefliger notes in \cite{Haefliger1970}  that the results of his paper may be generalized to the case where the pseudo-group \emph{$\Gamma_q$} is replaced by an arbitrary pseudo-group.

 In fact, the first general notion of geometric structure is usually attributed to Charles Ehresmann\index{Ehresmann, Charles} (1905-1979), who considered general geometric structures on manifolds defined by atlases of local coordinate charts with value in a homogeneous space. Ehresmann was Reeb's and Haefliger's doctoral advisor.

 In Thurston's setting, the first natural examples of $\mathcal{G}$-structures are the affine and projective structures: 
  If $G$ is the group of affine transformations of $\mathbb{R}^n$, then a $(G,\mathbb{R}^n)$-structure is called an \emph{affine structure}. A similar definition holds for projective structures. In his notes, Thurston gives examples of affine structures on the 2-torus, due to John Smillie, obtained by identifying opposite sides of quadrilaterals; see the latter's PhD dissertation \cite{Smillie1977}. It is also worth mentioning here Bill Goldman's Senior thesis, defended the same year \cite{Goldman-Senior}, in which he classifies  real projective structures on the 2-torus.  
 A few years after Thurston's notes were written, Sullivan and Thurston published the paper \emph{Manifolds with canonical coordinate charts: Some examples} \cite{Sullivan-Thurston} (1983). They write in the introduction: ``[\ldots] The geometry of the associated  developing maps\index{developing map} is a problem like the quantitative study of dynamical systems involving as it does the `infinite composition' of finitely many operations."
  In this paper,  Sullivan and Thurston study affine and projective structures in dimensions 2, 3 and 4.   The term ``manifold with canonical coordinate charts" that is used the title is another name for 
 a locally homogeneous geometric structure, that is, a manifold equipped with an atlas whose transition function belong to a certain pseudo-group of transformations that acts transitively on some model space. A particularly interesting example givan by Sullivan and Thurston is a projective structure on the 2-torus $\mathbb{T}^2$ whose developing map is not a covering and which gives
 an affine  structure on the 3-torus $\mathbb{T}^3$.

$G$-structures, as they first appear in the works of  Ehresmann,\index{Ehresmann, Charles}  were called\index{infinitesimal structure} \emph{infinitesimal structures}.\footnote{It seems that the name $G$-structure was first given by Chern, see \cite[p. 511]{Libermann} and the reference \cite{Chern}. Chern's paper builds on the works of S. Lie and \'E. Cartan. It is concerned with  complex-analytic  transformation groups  associated with  systems of partial differential equations.} Ehresmann\index{Ehresmann, Charles}  belonged to the first ``collaborators of Nicolas Bourbaki" (in fact, he was among the founders of the group). For Bourbaki, the notion of ``structure" was the basis of mathematics, and highlighting the notion of structure in various fields of mathematics was one of their basic undertakings.
Ehresmann's doctoral thesis, defended in 1934 and written under \'Elie Cartan, is titled \emph{Sur la topologie de certains espaces homogènes} (On the topology of certain homogeneous spaces), and he starts there by promoting the notion of geometric structure. In the Preface, we can read: ``The properties of a homogeneous space in which a transitive Lie group operates are simply an expression of the properties of this group. It would be interesting to know the relationship between the topology of such a space and the properties of its structure group. Our knowledge on this subject is still very incomplete. However, in his investigations of simple groups and symmetric homogeneous spaces, Mr. \'E. Cartan has come up with some remarkable results that reveal some of these relationships."\footnote{My translation from the French.} The first chapter of the thesis starts with the words: ``A homogeneous space of the most general type is a topological space whose homeomorphism group is transitive. In the following, we will reserve the name homogeneous space for an $n$-dimensional topological manifold that is transitively transformed by a finite continuous group $G$. This group will be called the structure group of the space and we will assume that it is a finite continuous Lie group." The idea of geometric structure is contained in this thesis. The thesis was published in the Annals of Math.  \cite{Ehresmann1934}. 
Other relevant papers by Ehresmann on geometric structures include  \cite{Ehresmann1952, Ehresmann1953, Ehresmann1954, Ehresmann1961}. 
In 1980, Saunders MacLane wrote a paper titled \emph{The genesis of mathematical structures, as exemplified
in the work of Charles Ehresmann} \cite{MacLane}. His paper\index{mathematical structure} has a philosophical tone; it starts with the words: ``Recently I have become interested in the need to revive the study
of the philosophy of mathematics." Then, later in the introduction: ``Thus we hope to examine how mathematical structures arise---not
principally how they actually happened to arise in the historical development, but how they logically and genetically must arise. [\ldots] The work of Charles Ehresmann provides splendid examples of the
genesis of mathematical structures."

Ehresmann  introduced the first definition of a topological or differentiable manifold in terms of  an atlas, see \cite{Ehresmann-atlas-1943}, 
and also \cite{Ehresmann-1947} which  contains the idea of a  local structure defined  by a groupoid of
allowable   transformations.
  Together with Georges Reeb, Ehresmann also introduced the notion of foliation, see \cite{Ehresmann-Reeb}. Among other geometric structures he introduced, let me mention the almost complex structures, the almost symplectic structures, the almost Hermitian structures, the almost K\"ahler structures and the almost quaternionic structures. See the report in \cite{Libermann}

Goldman published a captivating historical overview on geometric structures on manifolds, in particular flat affine, projective and conformal structures, see \cite{Goldman-Dani}.
He notes there that Ehresmann was heavily influenced by ideas of Lie, Klein, Poincaré and \'E. Cartan on this topics, and that he introduced the notions of
developing map, holonomy representation, normal structure and other related
notions which led to what Goldman calls the ``Thurston holonomy principle", or
the ``Ehresmann–Weil–Thurston holonomy principle". This principle establishes
a relation between the classification of geometric structures on a manifold and
the representation variety of its fundamental group into a Lie group, a topic of
extensive activity today. Goldman talks about ``Ehresmann's vision" which set
the context for Thurston's geometrization program for 3-manifolds. He mentions
the relation between the works of Poincaré, Lie, Klein and \'Elie Cartan and the
later developments in the theories of discrete subgroups of Lie groups, complex
projective structures, flat conformal manifolds and others. He shows how the notions of \emph{developing map} and \emph{holonomy representation} arose from the 
 idea of developing a $(G, X)$-manifold on its model space $X$. This is again a reformulation and  a generalization by Thurston of ideas of Ehresmann, whose origin lies in \'E. Cartan's work. In fact, the notion of geometric structure is close to that of flat Cartan connection and the idea of developing along a curve is closely related to that of parallel translation.  Goldman recalls 
that the question of classifying geometric structures on manifolds is a precise formulation of Klein's \emph{Erlangen program} (1872), which asks for a comparative study of the various known geometries at that time: Euclidean, spherical, hyperbolic, affine and projective. These geometries are among the first examples of geometric structures in the sense of Thurston.
 
 \section{Spheres and horospheres in hyperbolic space: Lobachevsky's insight} \label{s:spheres}
 
I would like to comment now on a sentence of Thurston concerning spheres and horospheres\index{horosphere} in hyperbolic geometry. He writes, in the chapter on geometric structures \cite[\S 3.8]{Thurston-Notes}: ``In the Poincaré
disk model, a hyperbolic sphere is a Euclidean sphere in the interior of the disk,
and a horosphere is a Euclidean sphere tangent to the unit sphere. The point $X$ of
tangency is the center of the horosphere."

My first comment is that Lobachevsky,\index{Lobachevsky, Nikolai Ivanovich}  without having at his disposal any Euclidean model of hyperbolic space, noticed that the geometric spheres in that space are spheres in the usual sense (they are endowed with a spherical geometry). My second comment is about horospheres:\index{horosphere} Lobachevsky showed that they carry a Euclidean geometry.   Let me explain this, and draw a consequence.

Lobachevsky studied geometric spheres in hyperbolic three-space, without models.
To investigate the geometry of such a sphere, he first defined the lines on this surface, by taking its intersections with (hyperbolic) planes passing through its center. He then defined angles between two lines to be the dihedral angles made by the planes that define them. The distance between two points on this sphere is the angle made by the geodesic rays starting at the the center of the sphere and passing through these points. This is analogous to the way one defines distances on a sphere embedded in Euclidean space.   Finally, Lobachevsky  defined triangles, as figures bounded by three    line segments which intersect pairwise at their endpoints. He then established the trigonometric formulae for these ``spherical" triangles and he noticed that these formulae 
coincide with the formulae of usual spherical geometry. Since in any geometry, the trigonometric formulae contain all the information on that geometry, Lobachevsky concluded that a geometric sphere in hyperbolic three-space carries the same geometry as the usual sphere in Euclidean space. He writes  (\cite{Pangeometry} p. 22 of the English translation):
``It follows that spherical trigonometry stays the same, whether we
adopt the hypothesis that the sum of the three angles of any rectilinear
triangle is equal to two right angles, or whether we adopt
the converse hypothesis, that is, that this sum is always less than
two right angles."

Regarding horospheres (which, in the \emph{Pangeomtry}\index{Lobachevsky, Nikolai Ivanovich!\emph{Pangeometry}}  Lobachevsky calls limit spheres),\index{limit sphere} he also considered triangles on such a surface: their sides are horocicrles (which he calls limit circles). He deduced that this 
 geometry  is Euclidean from the property that the angle sum of any
triangle is equal to two right angles. It is known indeed that this property
is equivalent to Euclid's parallel postulate. What is needed is to check that the
other axioms of Euclidean geometry are satisfied on the horosphere, and this can be done indeed.
Lobachevsky writes (p. 8 of the English translation \cite{Pangeometry}): ``It follows
from this theorem that the angle sum of any limit sphere triangle is equal to two right
angles and, consequently, everything that we prove in ordinary geometry concerning
the proportionality of edges of rectilinear triangles can be proved in the same manner in
the Pangeometry\index{Lobachevsky, Nikolai Ivanovich!\emph{Pangeometry}}  of limit sphere triangles."

Let me note now that Lobachevsky used these two facts, namely, the fact that in three-dimensional hyperbolic space the geometric spheres carry a spherical geometry and that the horospheres carry a Euclidean geometry, to prove the trigonometric formulae of hyperbolic geometry (again, in a model-free manner). The method is three-dimensional. It uses, for a given (right) triangle in  hyperbolic three-space, auxiliary triangles on a sphere and a horosphere,\index{horosphere} passing through a vertex, and the known trigonometric formulae for these triangles, see \cite[p. 22-30]{Pangeometry}. I have reviewed this in the article \cite{Papa-GB}.

All this is a hint to my belief that models are not an essential element of hyperbolic geometry. Everything can be done without models. I will talk again about this in \S \ref{s:Volume}, the section on the computation of volume.

\section{Polyhedra: Andreev's theorem}  \label{s:Andreev}

Subsection 13.6 of Thurston's notes is titled \emph{Andreev's theorem and generalizations}\index{Andreev's theorem} and it starts by the words: ``There is a remarkably clean statement, due to E. M. Andreev, describing hyperbolic
reflection groups whose fundamental domains are not tetrahedra."
Thurston states the following theorem, which he attributes to E. M. Andreev:

\bigskip

\noindent {\bf Theorem 13.6.1} (Andreev, 1967):

(a) Let $O$ be a Haken orbifold with $X_0=D^3$, $\Sigma_0=\partial D^3$. 

Then $O$ has a hyperbolic structure if and only if $O$ has no incompressible Euclidean suborbifolds.

(b) If $O$ is a Haken orbifold with $X_0=D^3$--(finitely many points) and $\Sigma_0= \partial X_0$, and if a neighborhood of each deleted point is the product of a Euclidean orbifold with an open interval (but $O$ itself is not such a product),
then $O$ has a complete hyperbolic structure with finite volume if and only
if each incompressible Euclidean suborbifold can be isotoped into one of the
product neighborhoods.

\bigskip
 
 There are two papers by Andreev related to Thurston's work, titled \emph{On convex polyhedra in Lobachevskii space} \cite{Andreev1},  and \emph{On convex polyhedra of finite volume in Lobachevskii space} \cite{Andreev2}.  In the first paper, Andreev  gives a classification of discrete isometry groups of hyperbolic space generated by reflections along hyperplanes with compact fundamental domain.  This is done by studying convex compact polyhedra in $n$-dimensional hyperbolic space whose dihedral angles are all non-obtuse. In the cases of Euclidean $n$-space and the $n$-sphere, the classification of the corresponding isometry groups was already known and is much easier, see Coxeter \cite{Coxeter1934}. Andreev starts by associating to any convex $n$-dimensional convex polyhedron $M$ an $(n-1)$-polyhedron which is combinatorially isomorphic to its boundary and which encodes its combinatorial type. He then proves several results, among them the fact that if two convex polyhedra $M$ and $M'$ have equivalent boundary polyhedra and if their corresponding dihedral angles are non-obtuse and equal, then $M$ and $M'$ are isometric. The last part of his paper is concerned with the applications to discrete group actions generated by reflections.

In the second paper, using the methods of the first paper, Andreev gives a complete classification of polyhedra of finite volume  with non-obtuse dihedral angles in three-dimensional hyperbolic space. (In the first paper, the polyhedra were assumed to be compact.) In particular, he obtains necessary and sufficient conditions for the existence of convex finite-volume polyhedra of a given combinatorial type in 3-dimensional hyperbolic space.

Andreev's characterization of compact hyperbolic polyhedra with all angles  non-obtuse 
is an essential ingredient of Thurston's Hyperbolization theorem for Haken 3-manifolds.

 The literature around Andreev's theorem became quite large after Thurston highlighted this work. 

A detailed proof of Andreev's result on the classification of three-dimensional compact hyperbolic polyhedra (other than tetrahedra) having non-obtuse dihedral angles was given in the PhD thesis of Roland Roeder \cite{Roeder}, which gave rise to the paper \cite{RHD} by 
Roeder, Hubbard and Dunbar. In the same paper, the authors correct a mistake in Andreev's existence proof. Other proofs of Andreev's theorem were given by several authors. We refer the reader to the bibliography in \cite{RHD}. Two of Thurston's students worked on generalizations, namely, Igor Rivin and Craig Hodgson. We mention in particular Rivin's thesis  \cite{Rivin}  in which he gave a characterization of compact, convex hyperbolic polyhedra, generalizing Andreev's result,  and the review of Rivin's result by Bowers in  \cite[\S 5.7.2]{Bowers}. See also 
Hodgson's paper \cite{Hod},  the review \cite{HRS} in the Bulletin of the AMS and the \emph{Erratum} \cite{HRS-E}.

Evgeny Mikhailovich Andreev  was a student of  \`Ernest Borisovich Vinberg. He obtained his PhD in 1970, at 
 Moscow State University 1970. The title of his doctoral dissertation is \emph{On convex polyhedra in Lobachevsky spaces}. It seems he has almost no published mathematical works after his two papers mentioned above.
 
 It is worth mentioning that there is a paper by Pogorelov on the realization of polyhedra in the hyperbolic 3-space as a right-angled polyhedron, which contains a particular case of Andreev's result and which seems to be poorly known, see \cite{Pogorelov}. Andreev, in his papers, did not quote Pogorelov.\footnote{I owe this remark to N. V. Abrosimov who, together with A. D. Mednykh, worked extensively on volumes of hyperbolic polyhedra.}

Milnor, in his paper on the first 150 years of hyperbolic geometry  \cite{Milnor-150}, notes that the importance of constructing discrete groups of hyperbolic isometries was already emphasized by Poincaré in his paper \cite{Poincare1882}, as a foundation for the study
of automorphic functions. Milnor notes that Poincaré referred to the existence of such examples due to Hermann Amandus Schwarz in \cite{Schwarz1873}. He also quotes Klein who, in his 1890 paper \cite{Klein1890}, inspired by examples due to Clifford, posed the following
 problem: ``\ldots to classify all connectivities which can possibly arise among closed
 manifolds of constant curvature." Here, ``connectivity" is understood to be a manifold, although no precise definition of an abstract manifold was given in the nineteenth century. 
 Milnor also indicates that Wilhelm Killing,  in his paper \cite{Killing1891}, pointed out the relation between Klein's classification problem and the study of discrete groups of isometries of Euclidean, hyperbolic or spherical  $n$-space. Killing called a group the quotient (that is, the manifold of constant curvature) of any  isometric, discrete and free action of a group on such a space gives  a \emph{Clifford-Klein space form}.

 My conclusion is that this work of Thurston on discrete isometry groups of hyperbolic space based on the study of polyhedra is again a continuation of a series of classical works.

\section{Volumes of polyhedra: Lobachevsky and Milnor}  \label{s:Volume}

Chapter 7, based on lectures by John Milnor, is titled  \emph{Computation of volume}, and is dedicated to volumes of hyperbolic ideal polyhedra.

 Let me start by recalling that Lobachevsky's memoirs on non-Euclidean geometry, since the first published one, contain extensive computations of hyperbolic volumes. 
 In fact, he developed the bases of a complete theory of differential and integral calculus in hyperbolic space (without models) which allowed him to compute areas and volumes of figures. 
 Computing the area or the volume of the
same figure in different ways allowed him to obtain equalities between definite integrals. Thus, computing volumes was a way of attracting attention on the geometry he invented, and showing that this geometry may be useful in the rest of the mathematical sciences. One may remember in this respect that the computation of definite integrals and finding  equalities between them was a very fashionable subject in those times. Cauchy in France and Ostrogradsky in Saint Petersburg,  were specialists in these matters.\footnote{In an ironic twist of fate, this was the reason why Lobachevsky's memoir \emph{On the elements of
geometry} was rejected, after he sent it in for approval by the  Russian Academy of Sciences, for publication in the Academy's journal.
 Rosenfeld, in his \emph{History of non-Euclidean geometry} \cite{Rosenfeld}, p. 208, quoting L. Kagan, reproduces a  
translation of Ostrogradsky's review, extracted from the St. Petersburg Academy
Archives:
``Having pointed out that of the two definite integrals Mr Lobachevsky
claims to have computed by means of his new method one is already known
and the other is false, Mr Ostrogradsky notes that, in addition, the work
has been carried out with so little care that most of it is incomprehensible.
He therefore is of the opinion that the paper of Mr Lobachevsky does not
merit the attention of the Academy."
 I have discussed the tragic fate of Lobachevsky in the book \cite{Pangeometry}. For this and for a general introduction to Lobachevsky's life and geometrical works, I refer the interested reader to that book.}
 The last section of Lobachevsky's \emph{Application
of imaginary geometry to certain integrals} \cite{Loba-Application} (1836)  contains a list of 50 
formulae for integrals.

E. B. Vinberg, in his paper \emph{Volumes of non-Euclidean polyhedra} \cite{Vinberg}, contains a detailed review of Lobachevsky's work on the computation of hyperbolic volume. The following facts are highlighted by Vinberg in his paper: Lobachevsky realized the importance of \emph{birectangular} hyperbolic tetrahedra, that is, tetrahedra $ABCD$ such that $AB$ is perpendicular to the plane $BCD$ and $DC$ is perpendicular to the plane $BCA$. In his memoir \cite{Loba-Application} he gave a remarkable formula that gives the volume of a birectangular tetrahedron\index{birectangular tetrahedron} in hyperbolic space in terms of the dihedral angles and the function which later became known   as the \emph{Lobachevsky function}.\index{Lobachevsky function} Now since any hyperbolic tetrahedron can be presented as an algebraic sum of birectangular ones, this formula is sufficient to give the volume of any hyperbolic tetrahedron. Vinberg writes that despite its importance, this formula ``was buried in the course of the noisy discussion about the foundations of geometry which developed in the second half of the 19th century and the beginning of the 20th." He adds that  the situation changed ``almost a hundred
and fifty years after Lobchevsky's works. This was connected with the works
of Thurston, who showed that to a significant extent the classification of
compact three-dimensional topological manifolds reduces to the classification
of compact hyperbolic spatial forms, that is, quotient manifolds of Lobachevsky
space under discrete groups of motions that act without fixed points and have
a compact fundamental domain."
Regarding Milnor's paper \cite{Milnor-150}, he writes:   ``His main observation was that the volume of a birectangular
tetrahedron with two vertices at infinity can easily be found by direct
integration in the Poincaré model. It turned out to be equal to half the
Lobachevsky function of one of the two equal acute dihedral angles of the
tetrahedron. (Of course, this also follows from Lobachevsky's formula.)
From this, in turn, it is easy to obtain a formula for the volume of any ideal
tetrahedron. The dihedral angles of such a tetrahedron at any pair of
opposite edges are equal, and its volume turns out to be equal to the sum of
the values of the Lobachevsky function\index{Lobachevsky function} at the dihedral angles of any vertex."

 It seems that the name \emph{Lobachevsky function} appears for the first time in Milnor's chapter of Thurston's notes (Chapter 7).
In this chapter,  Milnor introduces the Lobachevsky function,\index{Lobachevsky function} discussing its analytic properties, and he formulates two conjectures on number-theoretic properties of this function. He then shows how this function arises in computations of volumes of simplices and ideal simplices. In  Lobachevsky's memoir \emph{Application of imaginary geometry to certain integrals}  \cite{Loba-Application}, the function appears in a slightly different form. In the same chapter, Milnor discusses  examples of computation of volumes of hyperbolic manifolds that arise from knots and links  complements (essentially the figure-eight knot and the Whitehead link).  
 
 Thurston's motivations for studying three-dimensional hyperbolic volume were much more profound than Lobachevsky's. In his notes, based on Mostow’s rigidity theorem (of which he gave a new geometric proof),Thurston started by the fact that any complete finite-volume hyperbolic manifold is determined by its volume. Then, using techniques of hyperbolic Dehn surgery and of completion of incomplete structures, he obtained an infinite family of complete closed hyperbolic manifolds which are all non-homeomorphic. He showed that the set of volumes of hyperbolic three-manifolds is a closed non-discrete subset of the real line, and that it is well-ordered. He proved then that the volume function of hyperbolic three-manifolds is finite-to-one. He deduced from this that the set of volumes of hyperbolic three-manifold is indexed by countable ordinals and has the ordinal type of $\omega^{\omega}$. Furthermore, he showed that for any non-compact hyperbolic three-manifold of finite volume, there exists an infinite sequence of non-compact hyperbolic three-manifolds whose volume is strictly less than that of the original manifold, that approximate it in the sense of pointwise convergence of the representations of the fundamental groups in $\mathrm{PSL}(2, \mathbb{C})$, and whose volume converges to the volume of the original manifold. 
Thurston also addressed questions on number-theoretic properties of  hyperbolic volumes, some of which are still open. The work on hyperbolic volume which was motivated by Thurston's ideas is enormous and we shall not attempt to mention any part of the literature on this topic.

Now back to Chapter 7 of the notes.

Volume I of Milnor's \emph{Collected papers}, published in 1994, contains a chapter titled \emph{How to compute volume in hyperbolic space} \cite{Milnor-How-to-compute}.  Milnor indicates there that this chapter is a revised version of a manuscript he wrote in 1978, in relation with his notes included as Chapter 7 of Thurston's notes. He writes, in the introduction to this chapter: ``I was fascinated by Thurston's ideas, and corresponded with him regularly at this time. (We were both in Princeton, but seldom actually saw each other.)  The following paper is based on a handwritten manuscript, `Notes on hyperbolic volume' from 1978, which attempted to understand  how one can actually compute hyperbolic volume. Related material was circulated as part of the Thurston lecture notes."

In addition to the topics discussed in Chapter 7, Milnor's manuscript contains a discussion of groups generated by reflection in hyperbolic 3-space, related to Andreev's work, and a discussion of higher-dimensional computations.

Let me conclude with a curiosity concerning Lobachevsky. 
I recalled at the beginning of this section that he worked extensively on obtaining values of definite integrals, as an application of the geometry he discovered.
There is another topic on which he worked, as an application of his new geometry. Using the information on the stellar parallaxes that was available to him, he investigated the question of whether his ``theory of parallels holds or fails in nature". I refer the interested reader to the recent paper \cite{Berestovsky} by V. N.~Berestovski\u\i .

\bigskip 
 
\noindent {\bf Acknowledgement} The author is supported by the lnterdisciplinary Thematic lnstitute CREAA, as part of the ITl 2021-2028 program of the University of Strasbourg, CNRS and Inserm (ldEx Unistra ANR-10-IDEX-0002), and the French lnvestments for the Future Program.

\end{document}